\documentclass[11pt,laterpaper]{amsart}

\usepackage{amsfonts, amsmath, amssymb, amscd, amsthm, graphicx, verbatim}

\newtheorem{thm}{Theorem}[section]
\newtheorem{cor}[thm]{Corollary}

\newtheorem{quest}[thm]{Question}

\theoremstyle{definition}

\newtheorem{rem}[thm]{Remark}
\newtheorem{ex}[thm]{Example}

\newcommand{\e }{\varepsilon }
\def\co{\colon\thinspace}
\newcommand{\Aut}{\mathrm{Aut}}
\newcommand{\Gam}{\Gamma}
\newcommand{\Out}{\mathrm{Out}}

\newcommand{\Hom}{\mathrm{Hom}}

\newcommand{\FAb}{\mathrm{FAb}}

\renewcommand{\ll }{\langle\hspace{-.7mm}\langle }
\newcommand{\rr }{\rangle\hspace{-.7mm}\rangle }

\begin{document}

\title{Rips construction and Kazhdan property (T)}

\author{Igor Belegradek}
\address{Igor Belegradek, School of Mathematics, Georgia
Institute of Technology, Atlanta, GA 30332-0160}
\email{ib@math.gatech.edu}

\author{Denis Osin}
\address{Denis Osin, Department of Mathematics, The City College of
the City University of New York, 138th Street and Convent Ave., New
York, NY 10031} \email{denis.osin@gmail.com}
\thanks{The authors were supported by the NSF grants DMS-0503864 (Belegradek)
and DMS-0605093 (Osin).} \keywords{Hyperbolic group, Rips
construction, Kazhdan property (T)} \subjclass[2000]{20F67, 20F65,
20F28}
\date{}

\begin{abstract}
We show that for any non--elementary hyperbolic group $H$ and any
finitely presented group $Q$, there exists a short exact sequence
$1\to N\to G\to Q\to 1$, where $G$ is a hyperbolic group and $N$
is a quotient group of $H$. As an application we construct a
hyperbolic group that has the same $n$--dimensional complex
representations as a given finitely generated group, show
that adding relations of the form $x^n=1$ to a presentation of a
hyperbolic group may drastically change the group even in case
$n>> 1$, and prove that some properties (e.g. properties (T) and
FA) are not recursively recognizable in the class of hyperbolic
groups. A relatively hyperbolic version of this theorem is also
used to generalize results of Ollivier--Wise on outer
automorphism groups of Kazhdan groups.
\end{abstract}
\maketitle

\section{Introduction and main results}

Given any finitely presented group $Q$, the Rips construction
\cite{Rips} produces a short exact sequence
$$
1\to N\to G\to Q\to 1,
$$
where $G$ is a hyperbolic group and $N$ is finitely generated.
This result and its variations are a powerful source of
pathological examples in geometric group theory,
see~\cite{BMS,BriGru,Wis-res-fin,OllWis}. In this note the small
cancellation methods developed in \cite{SCT} and \cite{AMO} are
used to prove the following version of the Rips construction. (We
refer to the next section for definitions.)

\begin{thm}\label{rips}
Let $H$ be a non--elementary hyperbolic group, $Q$ a finitely
generated group, $S$ a subgroup of $Q$.  Suppose that $Q$ is
finitely presented with respect to $S$. Then there exists a short
exact sequence
\begin{equation}\label{seq}
1\to N\to G\stackrel{\e }{\to } Q\to 1,
\end{equation}
and an embedding $\iota \colon S\to G$ such that
\begin{enumerate}
\item[(a)] $N$ is isomorphic to a quotient group of $H$.

\item[(b)] $G$ is hyperbolic relative to the (proper) subgroup $\iota
(S)$.

\item[(c)] $\iota\circ \e \equiv id $.

\item[(d)] If $H$ and $Q$ are torsion free, then so is $G$.

\item[(e)] The canonical map $\phi\colon Q\to \Out (N)$ is injective and
$|\mathrm{Out} (N): \phi (Q)|<\infty $.
\end{enumerate}
\end{thm}

If $S=\{ 1\} $, we get the following corollary whose
applications are discussed in Section 3.

\begin{cor}\label{cor1}
For any non--elementary hyperbolic group $H$ and any finitely
presented group $Q$, there exists a short exact sequence
\rm(\ref{seq})\it, where $G$ is hyperbolic and $N$ is a quotient group of
$H$. Moreover, $G$ is torsion--free whenever $H$ is.
\end{cor}

\begin{rem}\label{rem1}
One can generalize Theorem~\ref{rips}
to the case where $H$ is a non--elementary group
hyperbolic relative to a collection of proper subgroups and $Q$ is
finitely generated and finitely presented relative to a collection
of subgroups $S_1, \ldots , S_m$. In this case the collection of
peripheral subgroups of $G$ will consist of isomorphic copies of
$S_1, \ldots , S_m$ and the images of the peripheral subgroups of
$H$ under the homomorphism $H\to N$. We leave the proof as an
exercise for the reader.

It is also worth noting that an alternative proof of Corollary
\ref{cor1} can be obtained by means of Ol'shanski{\u\i}'s
small cancellation theory over ordinary hyperbolic groups \cite{Ols}
without using any results of \cite{SCT} and \cite{AMO}.
\end{rem}

\section{Proof of Theorem \ref{rips}}

Recall that a group $G$ generated by a finite set $X$ is said to be
{\it finitely presented with respect to a subgroup $S$} if the
kernel of the natural map $\e\colon F(X)\ast S\to G$, where $F(X)$
is the free group with basis $X$, is the normal closure of a
finite subset $\mathcal R$ of $F(X)\ast S$. For instance, any group
is finitely presented relative to itself.

Given an element $w\in F(X)\ast S$, we denote by $|w|_{X\cup S}$ its
length with respect to the generating set $X\cup S$. A group $G$ is
{\it hyperbolic relative to $S$} if $G$ is finitely presented with
respect to $S$ and, in the above notation, there is a constant $C>0$
such that for any $w\in F(X)\ast S$ such that $\e (w)=1$, we have
\begin{equation}
w=\prod\limits_{i=1}^k f_i^{-1}R_i^{\pm 1}f_i, \label{prod}
\end{equation}
where $R_i\in \mathcal R$ and $f_i\in F(X)\ast S$ for
$i=1, \ldots, k$, and $k\le C|w|_{X\cup S}$.
In particular, this definition becomes the
definition of an ordinary hyperbolic group in the case $S=1$.

Let $G$ be a group hyperbolic relative to a subgroup $S$. Given a
subgroup $H\le G$, we denote by $H^0$ the set of all elements $g\in
H$ of infinite order such that that $g$ is  not conjugate to an
element of $S$. For every $g\in G^0$, there exists a unique maximal
virtually cyclic subgroup $E_G(g)$ of $G$ containing $g$
\cite{ESBG}. Moreover,
\begin{equation} \label{elem}
E_G(g)=\{ f\in G\; :\;
f^{-1}g^nf=g^{\pm n}\; {\rm for \; some\; } n\in \mathbb N\}.
\end{equation}

In this paper we use the notion of a suitable subgroup of a
relatively hyperbolic group introduced in \cite{SCT}. The most
convenient form of the definition is proposed in\cite{AMO}: A
subgroup $H\le G$ is called {\it suitable} if $H^0\ne \emptyset $
and
\begin{equation} \label{int}
\bigcap\limits_{g\in H^0} E_G(g)=\{ 1\}.
\end{equation}
In particular, (\ref{elem}) and (\ref{int}) imply that
every suitable subgroup $H\le G$ has trivial centralizer.

The key ingredient of the proof of Theorem \ref{rips} is the
following result obtained in \cite{SCT} (see Theorem 2.4 and its
proof there).

\begin{thm}\label{glue}
Let $G_0$ be a group hyperbolic relative to $S\le G_0$, $H$ a
suitable subgroup of $G_0$. Then for every finite subset $T=\{
t_1, \ldots , t_k\} \subset G_0$, there exist elements  $w_1,
\ldots , w_k\in H$ such that the quotient group
\begin{equation}\label{G}
G=G_0/\ll t_1w_1, \ldots , t_kw_k\rr ,
\end{equation}
satisfies the following conditions:
\begin{enumerate}
\item[\textup(i\textup)]  the restriction of the natural homomorphism $\xi\colon
G_0\to G $ to $S$ is injective;

\item[\textup(ii\textup)] $G$ is hyperbolic relative to $\xi (S)$;

\item[\textup(iii\textup)] if $G_0$ is torsion free, then $G$ is torsion free;

\item[\textup(iv\textup)] $\xi (H)$ is suitable in $G$;
\end{enumerate}
\end{thm}

\begin{proof}[Proof of Theorem \ref{rips}]
Suppose that $X\subset Q$ is a finite subset generating $Q$. Since
$Q$ is finitely generated and finitely presented with respect to
$S$, the group $S$ is also finitely generated by \cite[Theorem 1.1]{RHG}.
Let $Z$ denote a finite generating set of $S$. Assume also that the
kernel of the natural homomorphism $F(X)\ast S\to Q$ is the normal
closure of a finite subset $\mathcal R \subset F(X)\ast S$.
The free product $G_0=F(X)\ast S\ast H$ is hyperbolic
relative to $S$~\cite[Theorem 2.40]{RHG}.

It is well known  that every non--elementary hyperbolic group $H$
contains a unique maximal normal finite subgroup $K\le H$,
in fact $K$ is precisely the kernel of the $H$-action on the
boundary of $H$.
Thus passing to the quotient $H/K$ if necessary we may assume
that $H$ has no nontrivial finite normal subgroups. Since $H$ is
non--elementary, it is a suitable subgroup of $G_0$ and we may
apply Theorem \ref{glue} to the set
\begin{equation}\label{T}
T=\{ xyx^{-1},\ x^{-1}yx\; |\; x\in X\cup Z,\ y\in Y\}\cup \mathcal R,
\end{equation}
where $Y$ denotes a finite generating set of $H$. Let $G$ be the
corresponding quotient of $G_0$. By $N$ we denote the image of $H$
in $G$. Clearly (\ref{G}) and (\ref{T}) imply that $N$ is normal in
$G$ and $G/N\cong Q$. Assertions b)--d) of the theorem follow
immediately from Theorem \ref{glue}.

To prove (e) we first choose a non--elementary hyperbolic group $L$
with property (T). Ol'shanski{\u\i}'s results \cite{Ols} (see also
\cite{Ols-bl}) imply that there is a non--elementary hyperbolic
common quotient $H_1$ of $H$ and $L$. Thus replacing $H$ with $H_1$,
we may assume that $H$ has property (T), and hence so does $N$.
By part (iv) of Theorem \ref{glue}, $N$ is a suitable
subgroup, therefore (\ref{elem}) together with the definition of a
suitable subgroup implies that $N$ has trivial centralizer in $G$.
In particular, there is a canonical
injection $Q\to \mathrm{Out}(N)$. Now (e) follows as in
\cite[Section 5]{BelSzc}, and for completeness we sketch a proof
below.

Recall the definition of a relatively
hyperbolic group given in~\cite[Definition 1]{Bow}:
if $\mathcal G$ is a conjugacy-invariant family of infinite
finitely generated subgroups of a group $\Gamma$, then $\Gamma$ is
called {\it hyperbolic relative to} $\mathcal G$
if $\Gamma$ acts isometrically and properly discontinuously
on a complete locally compact hyperbolic metric space $\mathcal X$ such
that each point of the ideal boundary of $\mathcal X$ is either a conical
limit point or a bounded parabolic point of $\Gamma$, and
the elements of $\mathcal G$ are precisely the maximal
parabolic subgroups of $\Gamma$.

To continue the proof note that
if $S$ is infinite, then $G$ is hyperbolic relative
to $S$ in the sense of the above definition. Indeed, as we mentioned
before, since $G$ is finitely generated, so is
$S$~\cite[Theorem 1.1]{RHG}, and in this
situation the above definition of Bowditch is equivalent
to the definition introduced in~\cite{RHG}, which
is the definition we follow in this paper. The equivalence
of the definitions
is explained in appendix to~\cite{RHG}.
If $S$ is finite, then $G$ is hyperbolic~\cite[Corollary 2.41]{RHG},
so $G$ acts isometrically, properly discontinuously
and cocompactly on a locally compact complete hyperbolic space
$\mathcal X$.

Precomposing the inclusion $N\hookrightarrow G$ with
automorphisms of $N$, we get a sequence of injective homomorphisms
from $N$ to $G$, which defines a sequence of $N$-actions on $\mathcal X$.
These homomorphisms fall into finitely many
$G$-conjugacy classes, else a standard rescaling argument of
Bestvina-Paulin~\cite{Bes88, Pau88},
applied to the $N$-actions on $\mathcal X$,
produces a nontrivial action on an $\mathbb R$-tree
which contradicts the fact that $N$ has property (T). Since $N$ has
trivial centralizer in $G$, the canonical homomorphism
$G\to \mathrm{Aut}(N)$ in injective, and by the
previous sentence its image has finite index in $\mathrm{Aut}(N)$.
Let $G_\ast$ be the intersection of all the conjugates of (the image
of) $G$ in $\mathrm{Aut}(N)$. Note that $G_\ast$ is a
finite index normal subgroup of $\mathrm{Aut}(N)$ that contains $N$.
The kernel of the surjection of
$\mathrm{Out}(N)=\mathrm{Aut}(N)/N$ onto the finite
group $\mathrm{Aut}(N)/G_\ast$ is $G_\ast/N$.
Since $G_\ast/N\le G/N\le\mathrm{Out}(N)$, we conclude that
the image of $Q=G/N$ has finite index in $\mathrm{Out}(N)$.
\end{proof}

\section{Applications}

The idea behind many applications of the Rips
construction is to lift a (usually pathological)
property of $Q$ to $G$, and the extra control over $N$
provided by Theorem \ref{rips} makes the job easier.
\paragraph{\bf Prescribing linear representation.}
Here we build a hyperbolic group that has the same
set of $n$-dimensional complex representations as a given
finitely generated group; thus representation theory fails to
detect hyperbolicity.

\begin{thm}\label{thm: prescr linear}
For any finitely generated group $\Gam$ and any
integer $n>0$ there is a non-elementary
hyperbolic group $G$ and an epimorphism $G\to\Gam$
such that every representation $G\to GL_n(\mathbb C)$
factors as the composition of the epimorphism
$G\to\Gam$ and a representation $\Gam\to GL_n(\mathbb C)$.
\end{thm}
\begin{proof}
First, we build a finitely presented
group $Q$ and an epimorphism $Q\to \Gam$
such that any representation $Q\to GL_n(\mathbb C)$
factors through $Q\to \Gam$.
Let $\langle S|R\rangle$ be a presentation of $\Gam$
where $S=\{s_1,\dots,s_m\}$
and $R=\{R_1,\dots, R_k,\dots\}$.
The representation variety $\Hom (\Gam ,GL_n(\mathbb C))$
is the algebraic subvariety of the product of $m$ copies of
$GL_n(\mathbb C)$ defined by the relators $R_1,\dots, R_k,\dots$.
More precisely, $\Hom (\Gam ,GL_n(\mathbb C))$ is the set
of $m$-tuples of elements of $GL_n(\mathbb C)$ that are mapped
to the identity matrix by each $R_k$.
Let $V_k$ be the algebraic variety defined in the product of
$m$ copies of $GL_n(\mathbb C)$ by the first $k$ relators
$R_1,.., R_k$, and let $I_k$ be the
corresponding ideal (produced via Hilbert's Nullstellensatz).
Since polynomial rings over $\mathbb C$ are Noetherian,
the chain of ideals
$I_1,\dots,I_k,\dots $ stabilizes, i.e. $I=I_k$ for some $k$,
where $I$ is
the ideal that corresponds to $\Hom (\Gam ,GL_n(\mathbb C))$
(again via Hilbert's Nullstellensatz).
Hence the group $Q:=\langle s_1,...,s_m| R_1...R_k\rangle$
has the same representation variety as $\Gam$, i.e.
the inclusion $\Hom (\Gam ,GL_n(\mathbb C))\hookrightarrow
\Hom (Q ,GL_n(\mathbb C))$ induced by the epimorphism
$Q\to\Gam$ is a bijection.
In other words, any representation of $Q$ into
$GL_n(\mathbb C)$ factors through $Q\to\Gam$.

Next goal is to build a non--elementary hyperbolic
group with no nontrivial homomorphisms into $GL_n(\mathbb C)$.
It has been known for some time (see~\cite{Kap}) that
there exists a non--elementary hyperbolic group $K$
such that any homomorphism
$K\to GL_n(\mathbb C)$ has finite image.
In fact, by~\cite{Kap} any proper quotient of a cocompact
lattice in $Sp(r,1)$ with $r\ge 2$ has this property, and by small
cancellation theory there are proper quotients that are
non--elementary hyperbolic.
(Kapovich actually proves that
$K\to GL_n(\mathbb F)$ has finite image
for any field $\mathbb F$).
On the other hand,
Lubotzky observed in~\cite[Proposition 1.3]{Lub-superrig}
if a group $L$ has the property $\FAb$ (which means that
every finite index subgroup of $L$
has finite abelianization), then for each $n$ there is a finite
index subgroup
$L_n$ of $L$ such that for any representation
$\rho\co L_n\to GL_n(\mathbb C)$
the image $\rho(L_n)$ has connected closure.
Lattices in $Sp(r,1)$ have $\FAb$, because their finite index subgroups
have property (T).
Taking $L$ to be
any cocompact lattice in $Sp(r,1)$ with $r\ge 2$, we
get $L_n$ such that $\rho(L_n)$ has connected closure
for any $\rho\co L_n\to GL_n(\mathbb C)$.
Let $K_n$ be a non--elementary hyperbolic
proper quotient of $L_n$, so that
$\rho(K_n)$ is finite and has connected closure, and
therefore, is trivial.

Finally, if $G$ is the hyperbolic group produced by
the Rips construction
as in Corollary \ref{cor1} with input $H=K_n$
and $Q$ as in the first paragraph of the proof,
then any representation of $G$ into
$GL_n(\mathbb C)$ factors through
the epimorphism $G\to Q\to\Gam$.
\end{proof}

\begin{rem}
We do not know whether in Theorem~\ref{thm: prescr linear}
one can choose {\it the same} $G$ for all $n$.
The following representation-theoretic
properties of a group involve
its complex representations of {\it all} dimensions.
A finitely generated group $\Gamma$ is called
\begin{itemize}
\item
{\it representation rigid} if for every $n$
there are only finitely many pairwise non-isomorphic
$n$-dimensional irreducible representation
$\Gamma\to GL_n(\mathbb C)$.
\item
{\it representation superrigid} if there is a uniform upper bound
on the dimension of the Zariski closure of $\rho(\Gamma)$ for
all complex representations $\rho$ of $\Gamma$.
\end{itemize}
For brevity we just say {\it rigid} and
{\it superrigid}.
According to~\cite{BLMM} any superrigid group is
rigid. The Margulis superrigidity theorem implies
that irreducible lattices in semisimple groups of higher rank are
superrigid, and the same holds for lattices in
$Sp(n,1)$ and $F_4^{(-20)}$~\cite{Cor, Gro-Sch}.
Examples of non-linear rigid groups that are not superrigid can
be found in~\cite{BLMM}. However, we do not know any examples
of {\it finitely presented} rigid groups that are not superrigid.
(In the earlier version of this paper~\cite{Bel-Osin-arxiv},
we incorrectly assumed that
such examples were constructed in~\cite{BLMM}, and
we are grateful to Laszlo Pyber and David Fisher who independently
found the mistake).
As is explained in~\cite{Bel-Osin-arxiv}, if $Q$ is a
finitely presented rigid group that is not superrigid, and if
$K$ is a non--elementary hyperbolic proper quotient of a cocompact lattice
in $Sp(r,1)$ with $r\ge 2$, then the Rips construction
as in Corollary \ref{cor1} with input $H=K$
and $Q$ produces a non--elementary hyperbolic group $G$
that is rigid but not superrigid.
\end{rem}

\paragraph{\bf Prescribing minimal actions on Hadamard spaces.}
A $G$-action on a complete $CAT(0)$ space $X$ is called
{\it minimal} if $X$ is the only non-empty convex $G$-invariant
subspace. Let $\mathcal H$ be a class of complete $CAT(0)$ spaces
for which there exists a non-elementary hyperbolic
group $H$ with the property any isometric $H$-action on a space
$X$ in $\mathcal H$ has a non-empty fixed-point-set $X^H$.
Let $G$ be a hyperbolic group produced via the Rips construction
as in Corollary \ref{cor1} with input $H$, $Q$.

\begin{cor} \label{cor: cat0-minimal}
Any (isometric)
minimal $G$-action on every $X$ in $\mathcal H$
factors through the surjection $G\to Q$.
Thus $G$ and $Q$ have the same set of minimal actions on spaces
in $\mathcal H$.
\end{cor}
\begin{proof}
Since $N$ is a quotient of $H$, its fixed-point-set $X^N$
is non-empty, and also convex by~\cite[Corollary II.2.8]{BriHae}.
The subspace $X^N$ is stabilized by $G$ because
$N$ is normal in $G$, hence $X^N=X$ since the $G$-action is
minimal. Thus $N$ lies in the kernel of the $G$-action.
\end{proof}

\begin{rem}
The simplest example to which Corollary~\ref{cor: cat0-minimal}
applies is when $H$ is any Kazhdan hyperbolic group, and $\mathcal H$
is the class of $\mathbb R$-trees.
Note that the $G$-action is stable if
and only if the $Q$-action is stable.
In particular, Corollary~\ref{cor: cat0-minimal} implies that
studying stable actions of hyperbolic groups on $\mathbb R$-trees
is no easier than studying stable actions of finitely presented groups
on $\mathbb R$-trees.
Of course, Corollary~\ref{cor: cat0-minimal} also
has a relatively hyperbolic version, that shows in particular
that stable actions of finitely generated
relatively hyperbolic groups on $\mathbb R$-trees
can be as pathological as stable actions of finitely generated groups
on $\mathbb R$-trees.
\end{rem}

Here is a list of some other pairs
$(\mathcal H, H)$ to which Corollary~\ref{cor: cat0-minimal}
applies.

\begin{itemize}
\item
$H$ is any Kazhdan hyperbolic group, and $\mathcal H$
is the class of affine Hilbert spaces, and
cube (or more generally zonotopal)
$CAT(0)$ cell complexes~\cite{NibRee, HagPau};
\item
$\mathcal H$
is the class of Hadamard manifolds (i.e.
complete simply-connected Riemannian manifolds
of nonpositive sectional curvature), and
$H$ is the torsion-free Kazhdan hyperbolic group constructed
in~\cite{Gro-rand, IzeNay} such that any $H$-action
on a manifold in $\mathcal H$ has a fixed point;
\item
finite products of the spaces
listed above
provided the $H$-action preserves the product structure
(which is automatically true in most cases), and
fixes a point in every factor.
\end{itemize}

\paragraph{\bf Adding higher powered relations to hyperbolic groups.}
A well-known ``hyperbolic philosophy'' suggests that adding higher
powered relations to a hyperbolic group does not change the group
too much. For instance if $G$ is non--elementary hyperbolic and $g$
is an element of infinite order, than $G/\ll g^n\rr $ is also
non--elementary hyperbolic for all but finitely many $n$ \cite[Theorem 3]{Ols}.

Recall that a group is called {\it large} if it contains a finite
index subgroups admitting an epimorphism onto a non--abelian free
groups. In a recent paper~\cite{Lac}, Lackenby showed that for any
large group $G$ and any element $g\in G$, the quotient $G/\ll
g^n\rr $ is large for infinitely many $n$. Moreover, in case $G$
is free and non--abelian, the words ``infinitely many'' can be
replaced with ``all but finitely many''. (A short group theoretic
proof of these facts can be found in~\cite{OlsOsi}.) One may guess
that the same is true if $G$ is a large hyperbolic group and $g$
is an element of infinite order. (It is an easy exercise to
construct counterexamples in case we do not require $G$ to be
hyperbolic or allow the order of $g$ to be finite.) However, we have the
following.

\begin{cor}
There exists a large torsion--free hyperbolic group $G$ and an
element $g\in G$ such that $G/\ll g^n\rr $ is not large for all
odd $n\in \mathbb N$.
\end{cor}

\begin{proof}
Let us consider a short exact sequence (\ref{seq}), where
$$Q=\langle a,x,y\, |\, a^{-1}xa=y, \,
a^{-1}ya=x\rangle $$ and $N$ has property (T). Observe that the
subgroup $S$ generated by $ a^2, x, y$ has index $2$ in $Q$ and is
isomorphic to $\mathbb Z\times F_2 $, where $F_2$ is the free group
of rank $2$. Thus $Q$ is large and hence so is $G$.

Let $K=G/\ll g^n\rr $, where $g$ is a preimage of $a$ in $G$ and $n$
is odd. Suppose that a finite index subgroup $L$ of $K$ admits an
epimorphism $\e \colon L\to F$, where $F$ is a non--abelian free
group. Clearly $L$ splits into the short exact sequence $1\to U\to
L\to V\to 1$, where $U$ is a finite index subgroup of a quotient of
$N$ and $V$ is a subgroup of $Q/\ll a^n\rr\cong \mathbb Z/n\mathbb Z
\times \mathbb Z $. Since property (T) is inherited by quotients and
subgroups of finite index, $U$ has property (T). In particular $\e
(U)=1$ and hence the image of $L$ in $F$ is isomorphic to a quotient
group of the abelian group $Q/\ll a^n\rr $, which is a
contradiction.
\end{proof}

In~\cite[Lemma 12.2]{DT}, Dunfield and Thurston showed,
among other things, that the property to
have virtually positive first Betti number is inherited by all but
finitely many quotients of type $G/\ll g^n\rr $, $n\in \mathbb N$,
where $G$ is a certain amalgam of finite groups and $g$ is an
element of infinite order. Motivated by some long--standing
conjectures in 3--dimensional topology, they asked~\cite[Question 11.1]{DT}
if this is true
for arbitrary amalgam of finite groups and any element
of infinite order. Recall that every amalgam of finite groups
is hyperbolic. The next corollary shows that the
Dunfield--Thurston result can not be generalized to the case of
arbitrary hyperbolic groups.

\begin{cor}
There exists a torsion--free hyperbolic group $G$ and an element
$g\in G$ such that $\beta _1(G)>0$, but $G/\ll g^n\rr $ has property
(T) for all $n\in \mathbb N$. In particular, $G/\ll g^n\rr $
contains no subgroups of finite index with positive first Betti
number.
\end{cor}

\begin{proof}
It suffices to consider a short exact sequence (\ref{seq}), where
$Q$ is infinite cyclic and $N$ has property (T). If $g$ is a
preimage of a nontrivial element of $Q$ in $G$, then for any $n\ne
0$, $G/\ll g^n\rr $ is an extension of a quotient of $N$ by a finite
group. Hence $G/\ll g^n\rr $ has property (T) whenever $n\ne 0$.
Hence $G$ has no finite index subgroups with infinite
abelianization.
\end{proof}

\paragraph{\bf Recursively recognizable properties.}
A group-theoretic
property $P$ is called {\it recursively recognizable}
if there is an effective algorithm that decides from a finite
presentation of a group whether or not the group has the property $P$
(or more formally, $P$ is called {\it recursively recognizable}
if the set of finite presentations of groups with property $P$ is
recursive).

Recall that the isomorphism problem is decidable for torsion--free
hyperbolic groups \cite{Sel, GroDah}.
On the other hand many problems are known to be
undecidable for hyperbolic groups \cite{BMS,M}. For instance, the
property to be generated by a given number of elements is
recursively unrecognizable in the class of hyperbolic groups
\cite{BMS}. Here we provide some other examples. Let $\mathcal{P}$
be a property of groups. We consider the following conditions.

(1) $\mathcal{P}$ is not recursively recognizable in the class of
all finitely presented groups.

(2) $\mathcal{P}$ is preserved under taking quotients.

(3) If groups $A$ and $B$ have $\mathcal{P}$, then any extension
of $A$ by $B$ has $\mathcal{P}$.

(4) There exists a non--elementary hyperbolic group with property
$\mathcal{P}$.

\begin{ex}
Kazhdan's property (T) and Serre's property FA satisfy
(2)--(4). To show they are not recursively recognizable in the
class of all finitely presented  groups it suffices to note that
(T) implies FA \cite{Wat} and the free product $\mathbb Z/2\mathbb
Z \ast G$ has property (T) (or FA) if and only if $G\cong \{ 1\}
$. Since the property ``to be trivial'' is not recursively
recognizable in the class of all finitely presented groups, the properties
(T) and FA are not recursively
recognizable in the class of all finitely presented groups.
\end{ex}

\begin{cor}
If $\mathcal{P}$ is a property of groups satisfying (1)--(4), then
$\mathcal{P}$ is not recursively recognizable in the class of all
torsion--free hyperbolic groups.
\end{cor}

\begin{proof}
Let $H$ be a non--elementary hyperbolic group with property
$\mathcal{P}$. For every finitely presented group $Q$, we
construct a short exact sequence $1\to N\to G_Q\to Q\to 1$ as in
the theorem. By (2) and (3), the group $G_Q$ has $\mathcal{P}$ if
and only if $Q$ has $\mathcal{P}$. Since $G_Q$ is hyperbolic and
$\mathcal{P}$ is not recursively recognizable in the class of all
finitely presented groups, the result follows.
\end{proof}

\paragraph{\bf Kazhdan groups and their automorphisms.}
By a result of Paulin \cite{Pau}, the group
$\Out (H)$ is finite for every hyperbolic
group $H$ with Kazhdan property (T).
Answering a question of Paulin,
Ollivier-Wise~\cite{OllWis} constructed Kazhdan groups whose
outer automorphism groups are infinite. In fact, they showed that
any countable group $Q$ embeds into $\Out(N)$ for some Kazhdan
group $N$. In particular, $\Out (N)$ can fail to be finitely
presentable (in case $Q$ is not recursively presented).
As was noted in~\cite{BelSzc}, this result of Ollivier-Wise
implies that any
finitely presented group $Q$ can be embedded as a finite index
subgroup into $\Out(N)$ for some Kazhdan group $N$.

Our methods allow us to recover (and generalize) these results.

\begin{cor} \label{thm: out}
For any finitely generated group $Q$ there is split extension
$N\rtimes Q$ that is hyperbolic relative to $Q$, and such that
the canonical map $\phi\colon Q\to \Out (N)$ is injective,
$|\Out (N): \phi (Q)|<\infty $, and $N$ has Kazhdan property (T).
In particular, any countable group embeds into
$\Out(N)$ for some group $N$ with Kazhdan property (T).
\end{cor}

\begin{proof}
It suffices to apply Theorem \ref{rips} for a finitely generated
group $Q$, subgroup $S=Q$, and a non--elementary hyperbolic group
$H$ with property (T). The assertion "in particular" follows from
the fact that any countable group embeds into a finitely generated
one.
\end{proof}

\begin{quest}
If $H$ is a hyperbolic centerless group with property (T), then
$\Aut(H)$ is a finite extension of $H$ \rm\cite{Pau}\it, so
$\Aut(H)$ is also a hyperbolic group with property (T).
Hence the automorphism tower of $H$ consists of hyperbolic
groups with property (T). Does the tower terminate in
finitely many steps? This is true if $H$ is a centerless finite group
by a classical result of Wielandt \rm\cite[13.5.4]{Rob}.
\end{quest}

\section{Acknowledgements}
We are grateful to William Goldman for discussions on representation
varieties, and to the referee for editorial remarks.


\def\cprime{$'$} \def\cprime{$'$} \def\cprime{$'$}
\providecommand{\bysame}{\leavevmode\hbox to3em{\hrulefill}\thinspace}
\providecommand{\MR}{\relax\ifhmode\unskip\space\fi MR }
\providecommand{\MRhref}[2]{%
  \href{http://www.ams.org/mathscinet-getitem?mr=#1}{#2}
}
\providecommand{\href}[2]{#2}

\end{document}